\documentclass[10pt]{article}
\usepackage{fancybox}

\everymath{\displaystyle}
\usepackage{amsmath,amsfonts,amssymb}
\usepackage{amsthm}
\theoremstyle{plain}
\usepackage{graphicx,epsfig, color}
\usepackage{mathtools}
\usepackage{indentfirst}

\def\ss{\par\smallskip}
\def\ms{\par \medskip}
\def\bs{\par \bigskip}
\def\set#1{\{\,#1\,\}}

\def\inverse{^{-1}}

\def\bigvert{\ \big\vert\ }
\def\equaldef{ \buildrel\hbox{\scriptsize  def}\over= }

\def\and{\ \mbox{and} \ }
\def\be{\begin{enumerate}}
\def\ee{\end{enumerate}}
\def\bi{\begin{itemize}}
\def\ei{\end{itemize}}
\def\beqn{\begin{eqnarray*}}
\def\eeqn{\end{eqnarray*}}
\def\n{\mathbb N}

\def\F{\mathbb F}

\def\bc{\begin{center}}
\def\ec{\end{center}}
\def\ef{\end{frame}}
\def\ed{\end{document}}

\def\gF{\big(g(x), F(x)\big)}
\begin{document}

\newtheorem*{thm*}{Theorem}
\newtheorem{ex}{Example}[section]
\newtheorem{thm}[ex]{Theorem}
\newtheorem{defn}[ex]{Definition}
\newtheorem{cor}[ex]{Corollary}
\newtheorem{prop}[ex]{Proposition}
\newtheorem{lem}[ex]{Lemma}
\newtheorem{noat}[ex]{Note}
\newtheorem{rmk}[ex]{Remark}
\newtheorem{notation}[ex]{Notation}
\newtheorem*{qstn*}{Open question}
\noindent
\bc
{\bf\large Roots of Formal Power Series and
\ms
 New Theorems on Riordan Group Elements}
\ms
{\sc  Marshall M. Cohen}
\bs
{\small {\it Department of Mathematics,\\
 Morgan State University, Baltimore, MD}\\
Email: marshall.cohen@morgan.edu}
\ec
\hrulefill
\ms
\bc
\parbox{6in}{{\bf Abstract.}
{\sl Elements of the Riordan group $\cal R$ over a field $\F$ of characteristic zero are infinite lower triangular matrices which are defined in terms of pairs of formal power series.  We  wish to bring to the forefront, as a tool in the theory of Riordan groups, the use of multiplicative roots $a(x)^\frac{1}{n}$ of elements $a(x)$ in the ring of formal power series over $\F$ .  Using roots, we give  a Normal Form for non-constant formal power series, we prove a surprising simple Composition-Cancellation Theorem and apply this to show  that, for a major class of Riordan elements (i.e., for non-constant $g(x)$ and appropriate $F(x)$), only one of the two basic conditions for checking that $\gF$ has order $n$ in the group $\cal R$ actually needs to be checked.  Using all this, our main result is to generalize C. Marshall \cite{Marshall2017} and prove:  Given non-constant $g(x)$ 
satisfying necessary conditions,  there exists a unique $F(x)$, given by an explicit formula, such that $\big(g(x), F(x)\big)$ is an involution in $\cal R$.  Finally, as examples, we apply this theorem to ``aerated" series $h(x) = g(x^q),\ q\ \text{odd}$,  to find the unique $K(x)$ such that $\big(h(x), K(x)\big)$ is an involution. }  }
\ec
\hrulefill
\ms
\mbox{}\\
\small Key terms: Riordan group, formal power series, multiplicative roots of formal power series, involutions, group elements of order $n$.
\ms
\setlength{\parindent}{20pt}

\section {Introduction}
\noindent The Riordan group was introduced in \cite{SGWW1991} with applications to counting problems and combinatorial identities, and it has been of much interest in combinatorics. (See, for example, \cite{Barry2016}, \cite{Cohen2018b}, \cite{LMMS2012}, \cite{JL-Nkwanta2013}, \cite{SGWW1991}, \cite{Sprugnoli1994}.)
\subsection*{Riordan matrices}
\noindent Elements of the Riordan group $\cal R$  are infinite lower triangular matrices which are generalizations of Pascal's triangle. 
Let $\F$ be a field of characteristic zero. A {\it Riordan matrix } $A \in {\cal R} ={\cal R}(\F)$ is defined in terms of a pair of formal power series $g(x) \and F(x)$ and  is denoted as $A = \gF$ where
\beqn
g(x) & = & g_0 +g_1x +\cdots + g_nx^n + \cdots,\quad g_0\neq 0\\
F(x) & = & f_1x + f_2x^2 + \cdots + f_nx^n + \cdots, \quad f_1\neq 1\\
\eeqn%
If $A =\big(a_{m,\, n}\big)_{m, n \geq 0}$, then the zeroth column of $A$ has as generating function the formal power series $g(x)$ and, more generally,  the $n^{\rm th}$ column ($n \geq 0$) has as generating function the series $g(x)\cdot F(x)^n$. Thus,
\[A = \gF = \left[\begin{array}{ccccc}
       \big\vert &\big\vert&\big\vert&\cdots\\
       \big\vert &\big\vert&\big\vert&\cdots\\
       g(x) & g(x)\cdot F(x) & \ \ \ g(x)\cdot\big(F(x)\big)^2&\cdots\\
       \big\vert & \big\vert & \big\vert &\cdots\\
      \big\vert &\big\vert&\big\vert&\cdots\\
      \end{array}\right]\]
      and 
\[a_{m, n} = [x^m]\,\big(g(x)\cdot F(x)^n\big) \ \equaldef \ \text{coefficient of} \  x^m \ \text{in} \ g(x)\cdot F(x)^n.\]
\begin{ex}
{\bf \em Pascal's triangle} \ $P = \left[\,p_{n, k}\, \right]_{n, k \geq 0} = \left[{n \choose k}\right]_{n, k \geq 0} $ is given by
\[P = \left(\frac{1}{1 - x}, \frac{x}{1 - x}\right)\]
\end{ex}
\subsection*{The Riordan group}
\noindent Matrix multiplication of elements of ${\cal R}(\F)$ gives the rule
\[\gF\cdot\big(h(x), K(x)\big) = \Big(\,g(x)\cdot h\big(F(x)\big), K\big(F(x)\big)\Big).\]
This was introduced in \cite{SGWW1991} and is now called {\bf \em The Fundamental Theorem of Riordan Groups}. Under this multiplication $\cal R$ becomes a group with
\bi
\item identity element equal to $\big(1, x\big)$,
\item $\gF\inverse = \left(\frac{1}{g\big(\overline{F}(x)\big)}, \overline{F}(x)\,\right)$,
where $\overline{F}(x)$ is the compositional inverse of $F(x)$. (See Section 2.)
\ei
\subsection*{Involutions in $\cal R$}
\noindent An {\bf \em involution} in a group is an element of order two in that group.  Thus $\gF \neq (1, x)$ is an involution in $\cal R $
\begin{eqnarray}
& \iff & \gF^2 = \big(g(x)\cdot g\big(F(x)\big), \, F\big(F(x)\big)\big) = (1, x)\\
& \iff & \text{(a)}\  g(x)\cdot g\big(F(x)\big) = 1\ \and\ \text{(b)}\ F\big(F(x)\big) = x.
\end{eqnarray}
\noindent More generally, $\gF$ has {\it order $n$} in the group $\cal R$ iff $n$ is the least positive integer such that
\[\text{(a)}\  g(x) \cdot g\big(F(x)\big)\cdots \ \cdot g\big(F^{(n - 1)}(x)\big) = 1\ \and\ \text{(b)}\ F^{(n)}(x) = x,\]
where $F^{(k)}(x) \equaldef F\big(F( \cdots F(x)\big)\cdots\big)$,\ \  ($k$ times).

\subsection*{The goal of this paper}
\noindent In this paper we study the existence and uniqueness, given $g(x)$, of an $F(x)$ such that $\gF$ is an involution.  Thus we are vitally interested in the {\sl functional equations} (a) and (b) above, for power series. To study these, we highlight as a tool in the theory of Riordan groups the use of multiplicative roots $a(x)^\frac{1}{n}$ where $a(x) \in \F[[x]] = $ the ring of formal power series over $\F$ . 
\ms
\noindent In Section 2 we formulate basic facts concerning roots of formal power series, starting with \cite{Niven1969}, and we prove Theorem 2.5,  a surprising compositional cancellation theorem. 
\ms
\noindent In Section 3, we use this tool  to prove (Theorem 3.1) that in major situations only condition (a) need be checked in proving that $\gF$ has order $n$ and then to generalize generalize   C. Marshall \cite{Marshall2017} and prove (Theorem 3.3):  {\it Given a non-constant series
\[g(x) = g_0 + g_rx^r + g_{r + 1}x^{r + 1} + \cdots,\ \text{with}\  g_0 = \pm 1, g_r\neq 0 \and  r \ \text{odd},\]
then  there exists a unique $F(x) = - x + f_2x^2 + \cdots$, given by an explicit formula,  such that $(g(x), F(x))$ is an involution in $\cal R$. }
\ms
\section{Formal Power Series}
\subsection{Background and Notation}
\begin{notation} Given the field $\F$ of characteristic zero, we denote
\bi
\item $\F[[x]] = \set{g_0 + g_1x + g_2x^2 + \cdots \, \bigvert \, g_i \in 
\F}.$
\item $\F_0[[x]] = \set{g_0 + g_1x + g_2x^2 + \cdots \, \in \F[[x]] \bigvert \, g_0 \neq 0}.$
\item $\F_1[[x]] = \set{f_1x + f_2x^2 + f_3x^3 + \cdots \, \in \F[[x]] \bigvert \, f_1 \neq 0}.$
\item $\F_+[[x]] = \set{f_rx^n + f_{r + 1}x^{r + 1} + \cdots \, \in \F[[x]]  \bigvert \, r >0, f_r \neq 0}.$
\ei
\end{notation}
\noindent We start by summarizing well-known facts about the groups $\F_0[[x]]$,
$ \F_1[[x]]$.
(See, for example, \cite{Niven1969},\cite{Cohen2018a}, \cite{Gan2017} .) The Riordan group ${\cal R} = \cal{R}(\F)$ defined in the Introduction is  very rich algebraically because it is isomorphic to the semi-direct product of $\F_0[[x]] \and \F_1[[x]]$, where these are given very different group structures.
\be
\item $\F_0[[x]]$ is a group under multiplication: 
\ms
 If $g(x) = \sum_{n = 0}^\infty g_nx^n \in \F_0[[x] \and 
h(x) = \sum_{j = 0}^\infty h_jx^j \in \F_0[[x]],$ then
\[g(x)\cdot h(x) \equaldef \sum_{n = 0}^\infty\ \left(\sum_{i = 0}^n g_ih_{n -i}\right)x^n.\]
Moreover, if we write  $g(x) = g_0\cdot(1 + \frac{g_1}{g_0}x + \cdots) = g_0\cdot(1 + G(x))$, then
\beqn
\big(g(x)\big)\inverse & \equaldef & \frac{1}{g(x)} =
 \frac{1}{g_0}\Big(1 - G(x) + G(x)^2 - \cdots\Big)\\
 & = & \frac{1}{g_0} -  \frac{g_1}{g_0^2}x +\Big(-\frac{g_2}{g_0^2} + \frac{g_1^2}{g_0^3}\Big)x^2 + \cdots
\eeqn\\
\item $\F_1[[x]]$ is a group under substitution (= composition): 
\ms
If $F(x) = \sum_{n = 1}^\infty f_nx^n \in \F_1[[x]] \and 
K(x) = \sum_{j = 1}^\infty k_jx^j \in \F_1[[x]],$ then
\[(F\circ K)(x) \equaldef F\big(K(x)\big) \equiv f_1K(x) + f_2K(x)^2 + \cdots\]
The identity element is $\text{id}(x) = x$.
\ms
\noindent The compositional inverse of $F(x)$ is denoted $\overline{F}(x)$ and is given by
\[\overline{F}(x) = \frac{1}{f_1}x - \frac{f_2}{f_1^3}x^2 +\cdots \ \] 
(Solve $F\big(\overline{F}(x)\big)$ term-by-term or use the LaGrange Inversion Formula(\cite[p. 38]{Stanley1999}:\ 
\[[x^n]\,\overline{F}(x) = \frac{1}{n}[x^{n - 1}]\left(\frac{x}{F(x)}\right)^n.\ )\]
\ee
\subsection{Roots in $\F_0[[x]]$}
\noindent We state the foundational theorem on multiplicative roots of formal power series in Theorem 2.2 and then give  very useful Theorems ropositions on roots which follow from this.
\begin{thm} {\bf (Root Theorem 1)}\ \cite[\,Thm. 3, Thm. 17]{Niven1969}:
\ms
\noindent Suppose that $n\in \n$ and $a(x) = 1 + A(x),\ \ (A(x) \in F_+[[x]])$. Then
\be
\item [(a)] There exists a unique $b(x) = 1 + B(x),\ \ (B(x) \in F_+[[x]])$ such 
\ss
that
\ $\big(b(x)\big)^n = a(x)$.
\ms
\noindent We denote:\quad $b(x) \equiv \big(a(x)\big)^\frac{1}{n}$.
\item [(b)] More precisely, 
\ee%
\beqn
\text{If}\quad  a(x) & = &  1 + a_rx^r + a_{r + s}x^{r + s}  + a_{r + s +1}x^{r + s + 1} + \cdots ,\\
 & = & 1 + A(x)\ \ \quad\ \text{with}\  r, s > 0,\,  a_r\neq 0,\\
  then\  b(x) & = & \big(1 + A(x)\big)^{1/n}\\
    & = & \sum_{j = 0}^\infty {1/n \choose j}A(x)^j\ \text{(Generalized Binomial Theorem)}\\
\\
 \eeqn
 \beqn
 & = & 1 + \frac{1}{n}a_rx^r + b_tx^t  \ \text{where} \  t = \text{min}\, \{r + s, 2r\} \and\\
 && b_tx^t = 
\begin{cases} 
\frac{1}{n}a_{r + s}x^{r + s} & \text{if} \ r > s \\
&\\
\frac{1 - n}{2n^2}\cdot a_r^2\cdot x^{2r} & \text{if} \ r < s \\
&\\
\left(\frac{1}{n}a_{2r} + \frac{1 - n}{2n^2}\cdot a_r^2\,\right)x^{2r}  & \text{if} \ r = s \ \ 
\Box\\
\end{cases}\ \\
\eeqn 
\end{thm}

\begin{cor} {\bf (Normal form for non-constant $a(x)$)}
\ms
\noindent {\bf If} $a(x) \in \F[[x]]$ is non-constant, 
  \[a(x) = a_0 + a_rx^r + a_{r+1}x^{r + 1}  + \cdots, \  \ (r > 0, a_r\neq 0) \]
{\bf Then}  $a(x) $ may be {\bf uniquely} written as 
\[a(x) = a_0 + a_r\cdot \big(A(x)\big)^r\ \text{with $A(x)$ of the form} \ \ A(x) = x + \cdots .\]
{\bf Note: } $A(x) \in \F_1[[x]]$ with  $\overline{A}(x) = x + \cdots$.
\end{cor}
\noindent {\bf Proof:} $ a(x) = \ a_0 + a_r\cdot x^r \big(\underbrace{1 + \frac{a_{r+1}}{a_r}x + \frac{a_{r+2}}{a_r}x^2 + \cdots }_{\equiv \, \widehat{a}(x)}\big)$\\
\ms
\mbox{\hspace{.5in}} $  = a_0 +\ \  a_r\cdot x^r \widehat{a}(x) $\\
\ss
\mbox{\hspace{.5in}} $  = a_0 + a_r\cdot \Big(x\cdot\widehat{a}^{\ 1/r}(x)\Big)^r$\\
\ss
\mbox{ \hspace{.5in}} $= a_0 +\quad  a_r\cdot \big(A(x)\big)^r \ \text{where} \
 A(x) = x\cdot\widehat{a}^{\ 1/r}(x) \quad \Box$\\
\ms
\noindent The following theorem for roots of series $a(x) \in \F_+[[x]]$ (i.e., with constant term $= 0$) follows immediately from Theorem 2.2.
\begin{thm}\ {\bf (Root Theorem 2)}
\ms
\noindent {\bf If} $a(x) =  a_qx^q + a_{q+r}x^{q + r}  + \cdots, \  \ (q, r  > 0, a_q\neq 0)$\\
\mbox{\qquad \quad \ }$ = a_q\cdot\big(A(x)\big)^q\ \ \text{with}\ A(x) = x + \cdots $
\ms
\noindent {\bf then}
\be
\item  $B(x)$  
is a solution of the equation $B(x)^q = a(x)$ 
\[\iff \quad B(x) = b_1\cdot A(x),\ \ \text{where}\ b_1^q = a_q.\]
\noindent $ {\text{Indeed},} \qquad B(x) = b_1x + b_rx^r + \cdots\ \ \text{where}\ b_1^q = a_q$.
\item Given a $q^{\rm th}$ root $b_1$ of $a_q$, \ $B(x) = b_1\cdot A(x)$ is the  unique $q^{\rm th}$  root of $a(x)$ of the form $B(x) = b_1x + \cdots$.\hspace{1.2in} {$\Box$\ \ }
\ee
\end{thm}
\noindent {\bf Note:}  Our root $B(x)$ in the preceding theorem is an element of $\F_1[[x]]$ and has a compositional inverse $\overline{B}(x)$. It is this which often motivates us to take roots, as in the following Theorem.
\subsection{A Composition-Cancellation Theorem}
\begin{thm}  {\bf ( The Composition-Cancellation Theorem)}
\ms
\noindent Suppose that $g(x) = g_0 + g_rx^r + g_rx^{r + 1} + \cdots$ with $r \geq 1 \and g_r \neq 0$.
\ms
\noindent 
Let $a(x) = a_sx^s + a_{s + 1}x^{s + 1} + \cdots \and b(x) = b_sx^s + b_{s + 1}x^{s + 1} + \cdots,$
where $s \geq 1 \and a_s = b_s \neq 0$.  Then
\[g\big(a(x)\big) = g\big(b(x)\big)\ \implies \ a(x) = b(x).\]
\noindent
{\bf Note:} There is no assumption that $g(x)$ has a compositional inverse.
\end{thm}
\ms
\noindent {\bf Proof:}\  By Corollary 2.3, $g(x) = g_0 +g_r\cdot \Big(G(x)\Big)^r$ with $G(x) = x + \cdots \, \in \F_1[[x]]$. Thus
\beqn
 g\big(a(x)\big) & = & g\big(b(x)\big)\\
 \implies \quad \Big(G\big(a(x)\big)\Big)^r & =  & \Big(G\big(b(x)\big)\Big)^r  = (a_s)^rx^{rs} + \cdots\\
 & = & h_r\cdot\big(H(x)\big)^r\ \text{with}\\
 &&  \ \ h_r = (a_s)^r = (b_s)^r \and H(x) = x + h_2x^2 + \cdots\\
\eeqn
Thus, by Root Theorem 2 and the fact that $a_s = b_s$,
\beqn
G\big(a(x)\big) & = & a_sH(x) = b_sH(x) = G\big(b(x)\big) \\
\implies \quad 
\overline{G}\Big(G\big(a(x)\big)\Big) & = & \overline{G}\Big(G\big(b(x)\big)\Big)\\
\implies \hspace{.7in} a(x) & =&  b(x) \qquad \qquad \qquad \Box
\eeqn
\section {Theorems on Riordan Group Elements}
\subsection{Shortening the Proof that $\gF^n = (1, x)$}
\noindent From the definition of multiplication in the Riordan group $\cal R$ one immediately gets by induction the fact that, for $n\in \n$,  $\gF^n = (1, x) \ \text{(the identity matrix)}$ if and only if
\beqn
& \text{(a)} & \ g(x)\cdot g\big(F(x)\big) \cdots g\big(F^{(n - 1)}(x)\big) = 1\hspace{1.7in}\ \\
\text{and}\ & \text{(b)} & F^{(n)}(x) = x\\
\eeqn
\noindent In proving that an element of the group $\cal R$ has finite order, the proofs of (a) and (b) can each be quite intricate. (See \cite{Cohen2018a}, \cite{Cohen2018b}.) We use the Composition-Cancellation Theorem (2.5) to show that in very common situations, (b) follows automatically from (a).
\ms 
\begin{thm}  Suppose that $\gF \in {\cal R}$ with \\
$g(x)  =  g_0 + g_rx^r + g_{r + 1}x^{r + 1} + \cdots \and F(x)  = f_1x + f_2x^2 + \cdots$, 
where $g(x)$ is non-constant (\,$g_0 \neq 0, \, r > 0 \and g_r \neq 0$) and $(f_1)^n = 1$.
\ms
\noindent{Then:} \qquad $\gF^n = (1, x) \iff
g(x)\cdot g\big(F(x)\big)\cdots g\big(F^{(n - 1)}(x)\big) = 1.$
\end{thm}
\noindent {\bf Proof:} The necessity follows from the definition of multiplication in $\cal R$. To see sufficiency note that
\beqn
g(x)\cdot g\big(F(x)\big) \ \ \cdots \hspace{.3in} \cdots \hspace{.2in}  g\big(F^{(n - 1)}(x)\big) \hspace{.4in} \qquad & = & 1\\
\implies \qquad g\big(F(x)\big)\cdot g\big(F(F(x))\big) \cdots g\big(F^{(n - 1)}(x)\big)\cdot 
g\big(F^{(n)}(x)\big) &  = &1
\eeqn
\noindent But since  $ g\big(F^{(j)}(x)\big)$ has a multiplicative inverse and $F(x) = f_1x + \cdots$, \ we then have
\[g\big(F^{(n)}(x)\big)= g(x) \ \ \text{and}\ \ F^{(n)}(x)  = x + \cdots.\]                 
By The Composition-Cancellation Theorem (2.5), $F^{(n)}(x) = x.$ Thus
\[\gF^n = \big(g(x)\cdot g\big(F(x)\big)\cdots g\big(F^{(n - 1)}(x)\big) ,  F^{(n)}(x)= (1, x).\quad \Box\]

\subsection{The Unique Riordan Involution Determined by $g(x)$}
\noindent The method of using roots to solve functional equations, which is given in Corollary 2.3 and Theorem 2.5, is used in this section to determine, for given $g(x) \in \F_0[[x]]$,  exactly which (if any) $F(x)$ will make $\gF$ an involution  in the Riordan group. As a preliminary, the following Lemma gives necessary conditions for $\gF$ to be an  involution.
\ms
\begin{lem} Suppose that $\gF \in {\cal R}$ with 
\[ g(x)  =  g_0 + g_rx^r + g_{r + 1}x^{r + 1} + \cdots,\quad F(x)  = f_1x + f_2x^2 + \cdots.\]
Then 
\be
\item If $\gF$ is an involution then $g_0^2 = 1$.
\item If $g(x) = g_0$ then $\gF$ is an involution \\
$\iff \gF = (-1, x)$ \ \ or\ \  $\gF = \Big(\pm1, F(x) \Big)$,\\
 \mbox{\qquad\ } where $F(x) = - x + \cdots$ with $F(F(x)) =x$.
\item If $g(x) \neq g_0$ and $\gF$ is an involution ,\\
then  $r$ is odd,   $g_r \neq 0$ and $f_1 = -1$.
\ee
\end{lem}
\ms
\noindent {\bf Proof:} $\gF$ is an involution iff 
\bi
\item $x = F(F(x)) = f_1[f_1x + f_2x^2 + \cdots ] + f_2[f_1x + f_2x^2 + \cdots ]^2 + \cdots$\\
\mbox{\ \,} $= f_1^2x + \text{(higher powers)}$
\item $1 = g(x)\cdot g\big(F(x)\big) = g_0^2 + g_0g_r(1 + f_1^r)x^r + \cdots $.
\ei
Thus $f_1 = \pm1, \ g_0 = \pm1 \and \ \text{either}\  g_r = 0$ or $(r \ \text{is odd and} \ f_1 = -1)$.
The results of the Lemma follow. \quad $\Box$
\bs
\begin{thm} ({\bf Main Theorem, generalizing \cite{Marshall2017}})
\ms
\noindent Suppose that 
\beqn
g(x)  & =  & g_0 + g_rx^r + g_{r + 1}x^{r + 1} + \cdots,\ \text{with}\  g_0 = \pm1, r\ \text{odd},\ \and g_r \neq 0\\
& = & g_0 + g_r\cdot \big(G(x)\big)^r\ \text{with} \\
&& \qquad  G(x) = x + \cdots\ \text{(\,as in Corollary 2)}\\
&&\qquad \qquad \ = x\cdot\left(1 + \frac{g_{r + 1}}{g_r}x + \frac{g_{r + 2}}{g_r}x^2 + \cdots\right)^{1/r}.
\eeqn
Then there exists a unique $F(x) = f_1x + f_2x^2 + \cdots$ such that $\gF$ is an involution. This is  given by
\[F(x) =\overline{G}\left(\frac{-G(x)}{\big(g_0g(x)\big)^{1/r}}\right).\]
\end{thm}
\ms
\noindent {\bf Proof:} 
\beqn%
 && g(x)\cdot g\big(F(x)\big) = 1 \\
& \iff & \qquad \ \  g\big(F(x)\big) = \frac{1}{g(x)}\\
& \iff &   g_r\cdot \big(G(F(x))\big)^r = \frac{1}{g(x)} - g_0\\
& \iff & g_r\cdot \big(G(F(x))\big)^r = \frac{1 - g_0g(x)}{g(x)}\\
& \iff & g_r\cdot \big(G(F(x))\big)^r = \frac{-g_0g_r\cdot\big(G(x)\big)^r }{g(x)}\ \ (\text{since}\ g_0^2 = 1)\\
& \iff & \quad \big(G(F(x))\big)^r = \frac{-\big(G(x)\big)^r }{g_0g(x)}
\eeqn%
To solve the last functional equation, note that
\bi
\item $(-1)^r = -1$,
\item $\Big(G\big(F(x))\big)\Big)^r$\\
$ =\Big( [-x + f_2x^2 + \cdots] + G_2[-x + f_2x^2 + \cdots]^2 + \cdots\Big)^r$\\
\mbox{\hspace{.65in}}  $ = -\Big(x + \cdots\Big)^r$,
\item $g_0g(x) = 1 + g_0g_rx^r + \cdots$
\ei
\ms
 \noindent Thus $\gF \ \text{is an involution}$
\beqn &\iff  & g(x)\cdot g\big(F(x)\big) = 1 \hspace{.5in}\text{(by Theorem 3.1)} \\
&\iff &\ G\big(F(x)\big) = \frac{-\big(G(x)\big)}{\big(g_0g(x)\big)^{1/r}}\quad \ \text{(by Theorem 2.4)}\\
& \iff & F(x)  =\overline{G}\left(\frac{-G(x)}{\big(g_0g(x)\big)^{1/r}}\right). \quad \Box
\eeqn
\subsection{The Involution Determined by an Aerated $g(x)$} 
\ms
\noindent Perhaps the best known Riordan involution is 
$\gF = \left(\frac{1}{1 - x}, \frac{-x}{1 - x}\right)$ -- the Pascal matrix with negation of the odd-indexed columns.  If we ``aerate" $g(x) = {\textstyle \sum_{n = 0}^\infty x^n}$, adding blocks of zeroes to get \[h(x) =  \sum_{n = 0}^\infty x^{qn} = \frac{1}{1 - x^q},\ \  \, q \in \n,\ q\ \text{odd},\]
what series $K(x)$ makes $\big(h(x), K(x)\big)$ an involution? As an exercise in applying Theorem 3.3, one finds that
\begin{ex}
If $q \in \n$ is odd, the matrix 
\[\big(h(x), K(x)\big) = \left(\frac{1}{1 - x^q},\ 
\frac{-x}{\big(1 - x^q\big)^{1/q}}\right) =
\left(\sum_{n = 0}^\infty x^{qn},\ \sum_{n = 0}^\infty (-1)^{n + 1}{-\frac{1}{q} \choose n}x^{1 + qn}\right)\]
\ms
\noindent is an involution. 
\end{ex}
\noindent {\bf Proof:} Compute $h(x)\cdot h\big(K(x)\big) = 1$ and  apply Theorem 3.1. 
\quad $\Box$
\ms
\subsection*{The general case of an aeration}

\noindent More generally, suppose that $\gF$ is an involution and $q$ is a positive odd integer. Let $h(x) = g(x^q)$.  We use Theorem 3.3 to find $K(x)$ such that $\big(h(x), K(x)\big)$ is an involution.
For notational simplicity, we temporarily assume that $g(x) = g_0 + g_1x + g_2x^2 + \cdots$\ with\ $g_1 \neq 1$ ({\it i.e.,} that $r = 1$ in Theorem 3.3) and sketch how Theorem 3.3 leads us to $K(x)$. 
\ms
\noindent {\bf Notation (using Theorem 2.2):} If $a(x) = x^q + a_{q + 1}x^{q + 1} + \cdots,\ \ q\ \text{odd},$\\
then $\big(a(x)\big)^{1/q} \equaldef 
x\left(1 + a_{q + 1}x + a_{q + 1}x^2 + \cdots\right)^{1/q}$.
\ms
\noindent We have 
\bi
\item $g(x) = g_0 + g_1\cdot G(x)$ \ where $G(x) = x + \frac{g_2}{g_1}x^2 + \cdots$
\item $F(x) = \overline{G}\left(\frac{-G(x)}{g_0g(x)}\right)$ as in \cite{Marshall2017} or Theorem 3.3, since we are assuming $r = 1$.
\item $h(x) = g(x^q) = g_0 + g_1\cdot G(x^q) \equaldef h_0 + h_q\cdot \big(H(x)\big)^q$\\
\mbox{\qquad} where $h_0 = g_0,\ h_q = g_1 \and H(x) = \left[G(x^q)\right]^{1/q}$
\item $\overline{H}(x) = \big(\,\overline{G}(x^q)\,\big)^{1/q}$ since
\beqn
 H\Big(\big(\overline{G}(x^q)\big)^{1/q}\Big) & = & \left[G\left\{\Big(\big(\overline{G}(x^q)\big)^{1/q}\Big)^q\right\}\right]^{1/q}\\
 & = & \left[G\left\{\overline{G}(x^q)\right\}\right]^{1/q} = \left[x^q\right]^{1/q} = x\\
 \eeqn
 \item $K(x) = \overline{H}\left(\frac{-H(x)}{\big(h_0h(x)\big)^{1/q}}\right),\ \ \text{by Theorem 3.3}$  
 \ms
  $=  \left[\overline{G}\left(\frac{- G(x^q)}{(g_0g(x^q)}\right)\right]^{1/q}\  \text{(computation!),} 
 \quad =  \Big(F(x^q)\Big)^{1/q}$ (since $g_1 \neq 0$).
\ei
\ms
\noindent Having used Theorem 3.3 to discover the formula for $K(x)$ in the special case $r = 1$, the direct proof of the validity of this formula in general is quite simple, given its uniqueness:

\begin{thm} {\bf (The Aeration Theorem)}
\ms
\noindent Suppose that $\gF$ is an involution, where $g(x)$ is non-constant and  suppose that $h(x) = g(x^q)$ for some positive odd integer $q$.   Let $\big(h(x), K(x)\big)$ be the unique involution given by Theorem 3.3. \ 
 Then $K(x) =  \Big(F(x^q)\Big)^{1/q}$.
\end{thm}

\noindent {\bf Proof:}\ \ $h(x)\cdot h\left[\big(F(x^q)\big)^{1/q}\right]  =  
g(x^q)\cdot g\left[\Big(F(x^q)\Big)^{1/q}\right]^q = g(x^q)\cdot g\Big(F(x^q)\Big) = 1.$\\
\ms
\noindent The last term equals one by substituting $x^q$ into $g(x)\cdot g\big(F(x)\big) = 1$.  Since $K(x)$ is unique, it must equal $\Big(F(x^q)\Big)^{1/q}. \ \  \Box$

\end{document}